\documentclass[a4paper,12pt]{article}
\usepackage{epsfig}
\usepackage{color}
\usepackage{amsmath}
\usepackage{amsmath,amssymb}
\usepackage{comment}
\usepackage{ulem}
\usepackage{mathdots}

\usepackage{amsmath, amsthm, amssymb}
\usepackage{mathrsfs}
\usepackage{graphics}
\usepackage{epstopdf}
\usepackage{graphicx}
\usepackage{cite}
\usepackage{tikz}
\usepackage{pgfplots}
\pgfplotsset{compat=1.10}
\usetikzlibrary{shapes.geometric,arrows,fit,matrix,positioning}
\tikzset
{
    treenode/.style = {circle, draw=black, align=center, minimum size=1cm},
    subtree/.style  = {isosceles triangle, draw=black, align=center, minimum height=0.5cm, minimum width=1cm, shape border rotate=90, anchor=north}
}

\newtheorem{theorem}{Theorem}

\newlength{\boxedparwidth}  
{
\begin{center}
 \begin{tabular}{|@{\hspace{.15in}}c@{\hspace{.15in}}|}
 \hline 
   \begin{minipage}[t]{\boxedparwidth}
   \setlength{\parindent}{.25in}}%
   {\end{minipage} \\ \\ \hline 
 \end{tabular} 
\end{center} 
}

\title{Approximation of a damped Euler-Bernoulli beam model in the Loewner framework} 

\author{
I.V. Gosea \footnotemark[1]~ and  A.C. Antoulas \footnotemark[2] \ \footnotemark[1]
}

\date{}

\begin{document}

\maketitle

\renewcommand{\thefootnote}{\fnsymbol{footnote}}
%\footnotetext[1]{
%   Originally submitted September 2015.
%   This work was supported by NSF Grants CCF-1017401 and CCF-1320866, as well as 
%   DFG Grant AN-693/1-1.
%}
\footnotetext[1]{
	Data-Driven System Reduction and Identification Group, Max Planck Institute for Dynamics of Complex Technical Systems, Sandtorstrasse 1, 39106, Magdeburg, Germany
	({\tt gosea@mpi-magdeburg.mpg.de})
}
\footnotetext[2]{
	Department of Electrical and Computer Engineering, 
	Rice University, 6100 Main St, MS-366, Houston, TX 77005, USA 
	({\tt aca@rice.edu})
}

\begin{abstract} \smallskip
\noindent
The Loewner framework for model order reduction is applied to the class of  infinite-dimension systems. The transfer function of such systems is irrational (as opposed to linear systems, whose transfer function is rational) and can be expressed as an infinite series of rational functions. The main advantage of the method is the fact that reduced orders models are constructed using only input-output measurements. The procedure can be directly applied to the original transfer function or to the one obtained from the finite element discretization of the PDE. Significantly better results are obtained when using it directly, as it is presented in the experiments section.
\end{abstract}

\section{Introduction}

Model order reduction (MOR) is used to replace large, complex models of time dependent processes into smaller, simpler models that are still capable of representing accurately the behavior of the original process under a variety of conditions. The goal is an efficient and methodical strategy that yields lower dimensional systems which have (input-output) response characteristics close to those of the original system while requiring only a fraction of the large-scale simulation time and storage.

The motivation for MOR stems from the need for accurate modeling of physical phenomena that often leads to large-scale dynamical systems which require long simulation times and large data storage. For instance, one such example is that of discretization of partial differential equations over fine grids, which leads to large-scale systems of ordinary differential equations.

Throughout this work we exclusively consider interpolatory model reduction methods. These methods have initially emerged in numerical analysis and linear algebra and are related to rational interpolation. Roughly speaking, in the linear case, we will seek 
reduced models whose transfer functions match those of the original systems at selected frequencies, or interpolation points.

The milestone towards setting up the Loewner framework was Lagrange rational interpolation, i.e., constructing rational interpolants that use a Lagrange basis for the numerator and denominator polynomials. As shown in \cite{aci13}, this turns up to be a good idea mainly because the basis choice leads to algorithms that are numerically robust. Also, a rational function is more versatile than a simple polynomial because it has both poles (roots of the denominator) and zeros (roots of the numerator), and can model functions with singularities and highly oscillatory behavior more easily than polynomials can.

Using rational functions, models  that match (interpolate) given data sets of measurements are computed. In the context of linear dynamical systems, we start from data sets that contain frequency response measurements and we week the rational function that interpolates these measurements, i.e., we seek the linear dynamical system that models these measurements.

This is the approach covered by a class of widely used methods called Krylov or moment matching methods. More precisely, Krylov methods seek reduced-order models that interpolate the full-order transfer function and its derivatives (moments) $H(s_i) = \tilde{H}
(s_i), H^{(j)}(s_i) = \tilde{H} ^{(j)}(s_i) = 0$ at relevant frequencies $s_i$ for $j \geqslant 0$ . For these methods, interpolation typically leads to approximation of the whole transfer function over a wide interval of frequencies of interest.

The Loewner framework can be viewed as a step forward from the \textit{classical realization problem} which was covered extensively in the literature in previous years. The question arises as to whether such a problem can be solved if information about the transfer function at different points in the complex plane is provided. We will refer to this as the generalized realization problem. This problem can be stated as follows: given data obtained by sampling the transfer matrix of a linear system, construct a controllable and observable state space model of a system consistent with the data.

The classes of systems which can be treated in our approach range from linear to bilinear and also quadratic systems. For 
simplicity, consider SISO systems ($m,p=1$). We will deal only with linear systems whose dynamics are described in generalized state by:\
\vspace{-1mm}
\begin{equation}\label{eq:lin_sys}
\Sigma_L: \begin{cases} E\dot{x}(t) = Ax(t)+Bu(t), \\  \ \  y(t) = Cx(t). \end{cases} 
\end{equation}
We seek reduced systems of the form:
\vspace{-1.5mm}
\begin{equation}
\hat{\Sigma}_L: \begin{cases} \hat{E}\dot{\hat{x}}(t) = 
\hat{A} \hat{x}(t)+\hat{B} u(t), \\ \ \ \hat{y}(t) = \hat{C}\hat{x}(t), 
\end{cases}
\end{equation}
where $\hat{E} \in \mathbb{R}^{k \times k}, \ \hat{A} \in \mathbb{R}^{k \times k}, \ \hat{B} \in \mathbb{R}^{k \times m}, \ \hat{C} \in \mathbb{R}^{p \times k}, \ 
\hat{x}(t) \in \mathbb{R}^k$. The number of inputs and outputs $m$ and $p$, respectively, remain the same, while $k\ll n$.

Notice that the transfer function of a linear system is a rational function, or a rational matrix function in the case of multiple inputs and outputs, $m$ and $p$.

\section{The Loewner Framework}

The Loewner framework allows the identification of the underlying system directly from measurements. The key is that, in contrast to existing interpolatory approaches, larger amounts of data than necessary are collected and the essential underlying system structure 
is extracted appropriately. Thus, a basic advantage of this approach is that it is capable of providing the user with a trade-off between accuracy of fit and complexity of the model.  

The Loewner framework was developed quite recently and has its starting point in \cite{ajm07}. For further details we refer the reader to the thesis \cite{aci13} and to the tutorial  \cite{ali16}.

Given measurements composed of \textit{right interpolation data}: 
$\{(\lambda_i,r_i,w_i)|\lambda_i \in \mathbb{C},r_i \in \mathbb{C}^{m \times 1},w_i \in 
\mathbb{C}^{p \times 1}, i=\overline{1,\rho}\}$, \textit{and left interpolation data}: 
$\{(\mu_i,l_i,v_i)|\mu_i \in \mathbb{C},l_i \in \mathbb{C}^{1 \times p},v_i \in 
\mathbb{C}^{1 \times m}, i=\overline{1,\nu}\}$, one needs to construct a realization 
$\{E,A,B,C\}$, such that the \textit{right constraints}:\\ \vspace{-3mm}
\begin{equation*} \vspace{-1mm}
H(\lambda_i)r_i = [C(\lambda_iE-A)^{-1}B]r_i = w_i \ \forall \ i \in \ \{1,2,...,N\},
\end{equation*} %\vspace{-3mm}
and \textit{left constraints}: \vspace{-1mm}
\begin{equation*} \vspace{-1mm}
l_iH(\mu_i) = l_i[C(\mu_iE-A)^{-1}B] = v_i \ \forall \ i 
\in \ \{1,2,...,N\},
\end{equation*}
are satisfied.We define the Loewner matrix as follows:
\begin{eqnarray}
\mathbb{L} = \left[ \begin{array}{ccc}
\frac{v_1r_1- l_1w_1}{\mu_1-\lambda_1} & ... & \frac{v_1 r_\rho- l_1 
w_\rho}{\mu_1-\lambda_\rho} \\
\vdots & \ddots & \vdots\\
\frac{v_\nu r_1- l_\nu w_1}{\mu_\nu-\lambda_1} & ... & \frac{v_\nu r_\rho- l _\nu 
w_\rho}{\mu_\nu-\lambda_\rho}%
\end{array} \right] \in \mathbb{C}^{\nu \times \rho}.
\end{eqnarray}
\normalsize
The shifted Loewner matrix is defined as:
\begin{eqnarray}
\mathbb{L}_{\sigma} = \left[ \begin{array}{ccc}
\frac{\mu_1 v_1 r_1- \lambda_1 l_1 w_1}{\mu_1-\lambda_1} & ... & \frac{ \mu_1 v_1 r_\rho - 
\lambda_\rho l_1 w_\rho}{\mu_1-\lambda_\rho} \\
\vdots & \ddots & \vdots\\
\frac{\mu_\nu v_\nu r_1- \lambda_1 l_\nu w_1}{\mu_\nu-\lambda_1} & ... & \frac{\mu_\nu v_\nu r_\rho - \lambda_\rho l _\nu w_\rho}{\mu_\nu-\lambda_\rho}%
\end{array} \right] \in \mathbb{C}^{\nu \times \rho}.
\end{eqnarray}
\normalsize

It can be shown that the Loewner matrices defined above can be factored as follows: 
$\mathbb{L} = - \mathcal{O}E\mathcal{R}, \ \mathbb{L}_{\sigma} = - 
\mathcal{O}A\mathcal{R}$, $W = \mathcal{O}B$ and $V = C\mathcal{R}$ where $\mathcal{O}$ and $\mathcal{R}$ are \textit{tangential generalized reachability} and \textit{observability} matrices which are introduced in \cite{ajm07}. Based on this factorization, the next result arises:
\vspace{-3mm}
\begin{theorem}
If the Loewner pencil $(\mathbb{L},\mathbb{L}_{\sigma})$ is regular, then $\{\tilde{E}  = 
-\mathbb{L}, \ \tilde{A} = - \mathbb{L}_{\sigma}, \ \tilde{B}  = V, \ \tilde{C}  = W \ 
 \}$  is a realization of the data. Hence, we obtain that $H(z) = 
W(\mathbb{L}_{\sigma}-z\mathbb{L})^{-1}V$ is the required interpolant.
\vspace{-3mm}
\end{theorem}

Having the input data given by the matrices $\{\Lambda,M,V,W,L,R\}$, there exists an interpolant $H(z) = C(zE-A)^{-1}B$ if the 
following equality holds $\forall z \in \{\lambda_i\} \cup \{\mu_j\}$:
$$
\text{rank}(z\mathbb{L}-\mathbb{L}_{\sigma}) = \text{rank}[\mathbb{L} \ \  
\mathbb{L}_{\sigma}] = \text{rank}[\mathbb{L} \ \text{;} \ \mathbb{L}_{\sigma}] = k.
$$
\normalsize
If the Loewner pencil $(\mathbb{L},\mathbb{L}_{\sigma})$ is irregular, we compute the short SVD of the following matrices: 
$$
[\mathbb{L} \ \ \mathbb{L}_{\sigma}] = U 
\Sigma_l \tilde{Z}^{*}, \ \ \left[ \begin{array}{c} \mathbb{L} \\ \mathbb{L}_{\sigma} 
\end{array} \right] = \tilde{U} \Sigma_r Z^{*},
$$
where $U,\ Z \in \mathbb{R} ^{\rho \times k},\ 
\Sigma_l, \ \Sigma_r \in \mathbb{R} ^{l \times k}$. 

\begin{theorem}
The quadruple $(\tilde{E},\tilde{A},\tilde{B},\tilde{C})$, given by: $\{\tilde{E} = - U^{*} \mathbb{L} Z, \ \tilde{A} = -U^{*} \mathbb{L}_{\sigma} Z, \ \tilde{B} = U V, \ 
\tilde{C} = W Z\}$, is the realization of an approximate data interpolant.
\vspace{-1mm}
\end{theorem}
Thus, if we have more data than necessary, we can consider 
$(\mathbb{L},\mathbb{L}_{\sigma},V,W)$ as a singular model of the data.

\section{The Beam Model}

Model order reduction (MOR) plays a vital role in numerical
simulation of large-scale complex dynamical systems. These
dynamical systems are governed by ordinary differential equations
(ODEs), or partial differential equations (PDEs), or both.
To capture the essential information about the dynamics of the
systems, a fine semi-discretization of these governing equations
in the spatial domain is often required.

Since the solution of the PDE reflects the distribution of
a physical quantity such as the temperature of a rod or
the deflection of a beam, these systems are often called
distributed-parameter systems (DPS).

The transfer functions of DPS systems are irrational functions as opposed to the transfer functions of systems modeled by ordinary differential equations which are rational functions. Another difference is that the state space is infinite dimensional, usually
a Hilbert space. Consequently, DPS are also called
infinite-dimensional systems.  The analysis of rational and irrational transfer functions differ in a number of important aspects. The most obvious differences between rational and irrational transfer functions are the poles and zeros. Irrational transfer functions often have infinitely many poles and zeros.

The simplest example of transverse vibrations in a structure
is a beam, where the vibrations can be considered to occur only in one dimension.

Consider a homogeneous beam of length L experiencing small transverse vibrations. For small deflections the plane cross-sections of the beam remains planar during bending. Under this assumption, we obtain the classic Euler-Bernoulli beam model for the deflection $w(x, t)$:
\begin{equation}
\frac{\partial^2 w(x,t)}{\partial t^2} + EI \frac{\partial^4 w(x,t)}{\partial x^4} = 0.
\end{equation}

Consider the Kelvin-Voigt damping model, sometimes referred to as the Rayleigh damping which leads to the PDE: 
\begin{equation}
\frac{\partial^2 w(x,t)}{\partial t^2} + EI \frac{\partial^4 w(x,t)}{\partial x^4}+ c_d I  \frac{\partial^5 w(x,t)}{\partial x^4 \partial t} = 0,
\end{equation}
with the boundary conditions at the clamped end:
\begin{equation}
w(0,t) = 0, \ \ \frac{\partial w(0,t)}{\partial x} = 0,
\end{equation}
and, additionally, at the free end:
\begin{equation}
EI \frac{\partial^2 w(L,t)}{\partial x^2} + c_d I \frac{\partial^3 w(L,t)}{\partial x^2 \partial t} = 0, \ \ -EI \frac{\partial^3 w(x,t)}{\partial x^3} - c_d I \frac{\partial^4 w(L,t)}{\partial x^3 \partial t}  = u(t),
\end{equation}
where $u(t)$ represents an applied force at the tip. Also take the observation (output) to be $y(t) = \frac{\partial w}{\partial t}(L,t)$and let $z(t) = w(L,t)$. 

This model corresponds to the case of a clamped-free beam with shear force control covered in \cite{cm09}. By taking different boundary conditions, one might analyze different models such as clamped-free beam with torque control or pinned-free beam with shear force control. For this study, we restrict our attention only to the first model. As derived in \cite{cm09}, the transfer function can be written as follows:
\begin{equation}\label{tref}
H_{orig}(s) = \frac{s N(s)}{(EI+sc_dI)m^3(s)D(s)},
\end{equation}
in terms of the following nonlinear functions:
\vspace{-2mm}
\begin{eqnarray*}
m(s) &=& \Big{(} \frac{-s^2}{EI+sc_dI} \Big{)} ^{\frac{1}{4}}, \\
N(s) &=& \cosh(Lm(s))\sin(Lm(s))- \sinh(Lm(s))\cos(Lm(s)), \\
D(s) &=& 1 + \cosh(Lm(s))\cos(Lm(s)).
\end{eqnarray*}
The exact derivation of these types of irrational transfer functions is described in \cite{book2}.\\
\noindent 
The poles of the transfer function (\ref{tref}) are the solution of the equation: $s^2+c_dI \alpha_k^4s+EI\alpha_k^4 = 0$ and can be written as:
\begin{equation}
\mu_{\pm k} = \frac{-c_dI \alpha_k^4 \pm \sqrt{(c_dI)^2\alpha_k^8-4EI\alpha_k^4}}{2},
\end{equation}
where the coefficients $\alpha_k$ are the real positive roots of the hyperbolic equation in $\alpha$: $1+\cosh(L\alpha) \cos(L\alpha) = 0$. When $k \rightarrow \infty$, these values converge to $\frac{(2k+1)\pi}{2L}$.

The complex poles $\mu_{\pm k}$ approach the imaginary axis as $c_d \rightarrow 0$ (the case for 0 damping). Notice that there is also a real pole at $-\frac{E}{c_d}$.

The zeros of the transfer function are the solution of the equation: $s^2+c_dI \gamma_k^4s+EI\gamma_k^4 = 0$ where $\gamma_k$'s are the roots of the the equation in $\gamma$: $\cosh(L\gamma)\sin(L\gamma)-\sinh(L\gamma)\cos(L\gamma) = 0$. When $k \rightarrow \infty$, these values converge to $\frac{(4k+1)\pi}{4L}$.

The original transfer function is irrational and is written as an infinite partial fraction expansion:
\begin{equation}
H_{orig}(s) = \displaystyle \sum_{k=1}^{\infty} \frac{r_k}{s-\mu_k}+ \frac{r_{-k}}{s-\mu_{-k}} = \displaystyle \sum_{k=1}^{\infty} \frac{4s}{s^2+c_dI\alpha_k^4s+EI\alpha_k^4},
\end{equation}
where $r_k = \frac{4 \mu_k}{\mu_{-k}-\mu_k}, r_{-k} = \frac{4 \mu_{-k}}{\mu_{-k}-\mu_k} $ are the residues of the poles $\mu_{\pm k}$.\\

We consider that the beam has length $L = 0.7$m and that it is build out of aluminum. Hence the Young modulus elasticity constant is taken to be  $E = 69 GPa = 6.9 \times 10^{10} \frac{N}{m^2}$. Then, the height and base of the rectangular cross section of the beam are taken to be $h = 8.5$mm and $b = 70$mm. Then, since the moment of inertia can be calculated precisely in terms of these two quantities: $I = \frac{bh^3}{12}$, it follows that: $I = 3.58 \times 10^{-9} m^4$. Finally, take the damping coefficient to be $c_d = 5 \times 10^{-4} \frac{Ns}{m^2}$.

One way to proceed in modeling the dynamics of the beam is by computing a finite element discretization (get rid of the variable x by approximating the spatial derivatives).  Assume that the beam has been divided in $N+2$ intervals of length $h = \frac{L}{N+2}$; the resulting variables are $w(kh,t), \ k \in \{0,1,...,N+2\}$. Consider $w_k = 0$ for $k \geqslant N+3$. \\
\noindent
For simplicity, instead of $w(kh,t)$, we are using $w_k$. Approximate the time derivatives of $w(x,t)$ at $x = kh$ with the following divided difference formulas:
$$
w^{(1)}_k \approx \frac{w_{k+1}-w_{k}}{h}, \ \  w^{(2)}_k \approx \frac{w_{k+1}-2w_{k}+w_{k-1}}{h^2},
$$
$$
 w^{(3)}_k \approx \frac{w_{k+2}-3w_{k+1}+3w_{k}-w_{k-1}}{h^3}, \  w^{(4)}_k \approx \frac{w_{k+2}-4w_{k+1}+6w_{k}-4w_{k-1}+w_{k-2}}{h^4}.
$$
The first two boundary conditions give us the following constraints: $w_0 = 0$ and $\frac{w_1-w_0}{h} = 0$ which means that we can eliminate the first two variables from the state vector ($w_0 = w_1 = 0$).\\
\noindent
The other two boundary conditions give us:
\begin{align*}
 EI \Big{(} \frac{w_{N+2}-2w_{N+1}+w_{N}}{h^2} \Big{)}  +  c_d I \Big{(} \frac{\dot{w}_{N+2}-2\dot{w}_{N+1}+\dot{w}_{N}}{h^2} \Big{)} &= 0, \\
 -EI \Big{(} \frac{w_{N+2}-3w_{N+1}+3w_{N}-w_{N-1}}{h^3} \Big{)} - c_d I \Big{(} \frac{\dot{w}_{N+2}-3\dot{w}_{N+1}+3\dot{w}_{N}-\dot{w}_{N-1}}{h^3} \Big{)} &= u.
\end{align*} 
It follows that:
\begin{align*}
 EIw_{N+1}+ c_dI\dot{w}_{N+1} &= EI(2w_N-w_{N-1})+c_dI(2\dot{w}_N-\dot{w}_{N-1})+ h^3 u, \\
 EIw_{N+2}+ c_dI\dot{w}_{N+2} &= EI(3w_N-2w_{N-1})+c_dI(3\dot{w}_N-2\dot{w}_{N-1})+ 2h^3 u.
\end{align*}
Hence, we are left with $N-1$ variables in the state vector: $v = [w_2,w_3,...,w_{N}]^T \in \mathbf{R}^{N-1}$. For each variable, write down the corresponding ODE: $w_k^{(2)}+EI w_k^{(4)}+c_dI \dot{w}_k^{(4)} = 0, \ 2 \leqslant k \leqslant N$ where $w_k^{(j)}$ is defined above. Then, the next step is to rewrite these ODE's into matrix format: $M \ddot{v} + U \dot{v}+ K v = f u$. Since $y(t) = w (L,t)$, the output is chosen as follows: $y = \dot{w}_{N+2}$.

\section{Experiments}

In the past years, numerical approximation of large-scale dynamical systems has been increasingly covered in the literature. Methods such as \textit{Balanced Truncation}, \textit{Rational Krylov}, \textit{Optimal} $\mathcal{H}_2$ \textit{approximation} and \textit{Modal Truncation} are just some of the contributions. We refer the reader to \cite{book1} as an extensive survey on the above mentioned methods. In this short study we will compare the results of the Loewner method against the results obtained when applying modal truncation only.

Using the Loewner framework, we will construct reduced order linear models approximating the original transfer function of the clamped-free beam with shear force control model in section 3. We will do so directly by using measurements of the original transfer function or by first performing a finite element semi-discretization of the PDE and then taking measurements. The dimension of all reduced models is taken to be constant throughout all experiments: $r=32$. For the choice of the parameters in Section 3, proceed by sampling $H_{orig}$ on a logarithmic frequency range between  $[10^0,10^7]j$. The frequency response of the original transfer function is depicted in Fig.\;1.

\begin{figure}[h] \label{fig2} \vspace{-2mm}
	\begin{center}
		\includegraphics[scale=0.3]{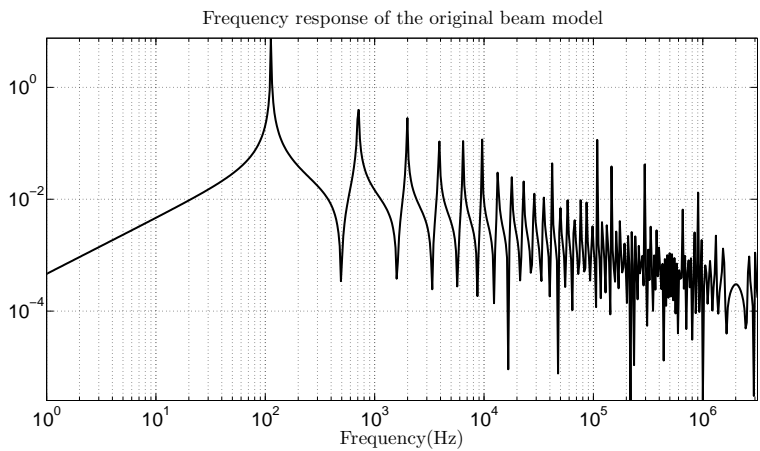}
		\vspace{-3mm}
		\caption{Frequency response of the beam model}
	\end{center} \vspace{-6mm}
\end{figure}

%Only when reaching $10^{10}$ Hz the peaks seem to stagnate at a value below $10^{-5}$ as it can be observed in the next plot - we increase the frequency range to $[10^0,10^{12}]j$.
% 
%\begin{center}
%\includegraphics[scale=0.3]{plot_beam2.eps}
%\end{center}

The first 32 poles and the first 31 zeros of the original transfer function are intertwined as we can observe in Fig.\;2. Additionally, notice that the absolute value of the imaginary part is less than $2.5 \times 10^{-5}$ for all $32$ poles - and it is slowly increasing when taking less and less dominant poles. The imaginary part of the dominant pair of poles is $-4.6113 \times 10^{-11}$. Notice that as the damping constant $c_d$ becomes smaller and smaller, the poles approach more and more the $j\omega$ axis. 

\begin{figure}[h] \label{fig2} \vspace{-3mm}
	\begin{center}
		\includegraphics[scale=0.32]{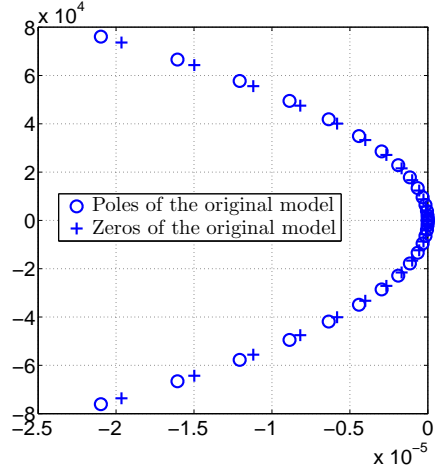}
		\vspace{-3mm}
		\caption{Poles and zeros distribution of the beam model}
	\end{center} \vspace{-6mm}
\end{figure}

We are going to use our data driven MOR technique, i.e. the Loewner method, to find a linear model that interpolates the original transfer function at some chosen frequencies in some particular 'interest' range. In the next experiments, this range is fixed to be $[10^1,10^{4.5}]$. The number of interpolation points is chosen to be 400. The plot in Fig.\;3 depicts the samples taken in the aforementioned interval. The linear model that is constructed based on these input-output measurements via the Loewner method is going to be singular (from a numerical point of view).
 
\begin{figure}[h] \label{fig3} \vspace{-3mm}
	\begin{center}
		\includegraphics[scale=0.3]{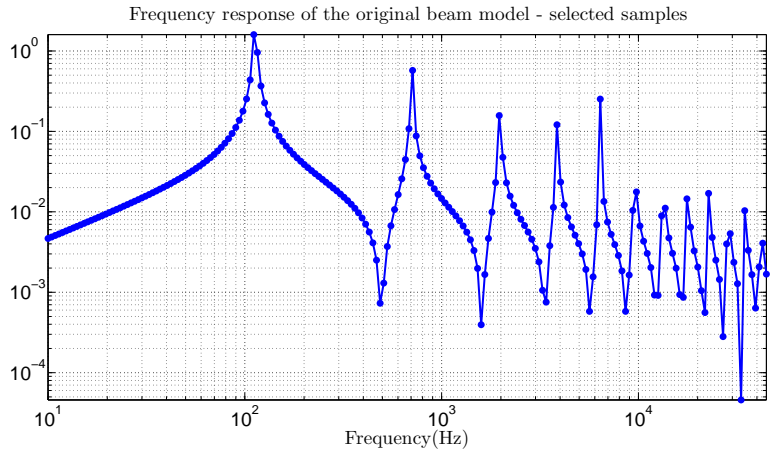}
		\vspace{-3mm}
		\caption{Samples corresponding to the frequency response of the beam model}
	\end{center} \vspace{-6mm}
\end{figure}

 The decay of the singular values of the Loewner matrix (see Fig.\;4) represents a good indicator that we can compress the model to a much smaller dimension (i.e. the $40^{\rm th}$ singular value is below machine precision $10^{-16}$).

\begin{figure}[h] \label{fig4} \vspace{-3mm}
	\begin{center}
		\includegraphics[scale=0.3]{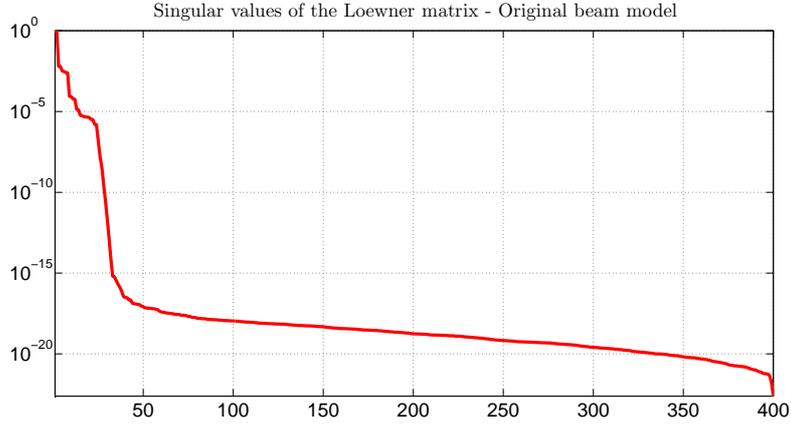}
		\vspace{-3mm}
		\caption{Decay of the Loewner matrix singular values}
	\end{center} \vspace{-6mm}
\end{figure}

We choose reduction order 32. Note that the initial samples are accurately matched.  Additionally, outside the interest range, the fitted linear model is not able to faithfully reproduce the peaks.

\begin{figure}[h] \label{fig5} \vspace{-3mm}
	\begin{center}
		\includegraphics[scale=0.3]{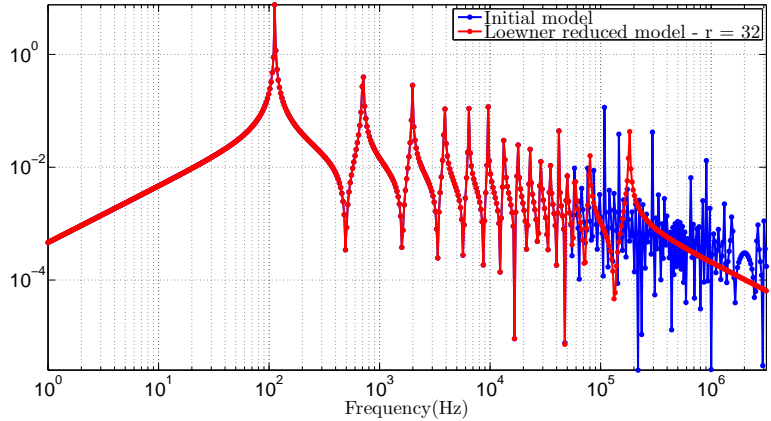}
		\vspace{-3mm}
		\caption{Original and Loewner model - Frequency response comparison}
	\end{center} \vspace{-6mm}
\end{figure}

Notice that the zeros and poles of our reduced linear model are matching the ones of the original transfer function. Only the least three dominant poles/zeros are not accurately reproduced.

\begin{figure}[h] \label{fig6} \vspace{-3mm}
	\begin{center}
		\includegraphics[scale=0.3]{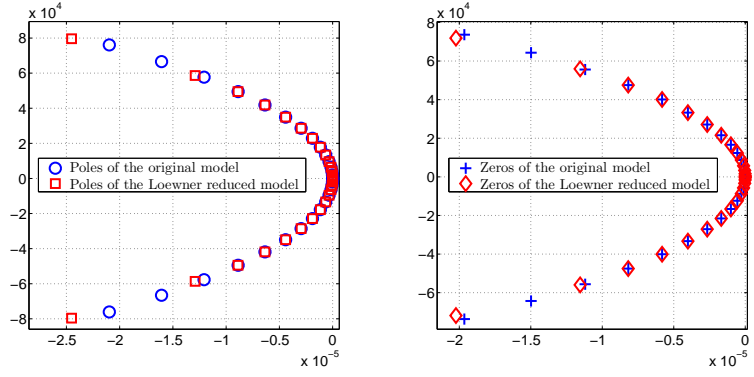}
		\vspace{-3mm}
		\caption{Original and Loewner model - poles and zeros comparison}
	\end{center} \vspace{-6mm}
\end{figure}

We want to compare our reduction method with other available methods for reducing systems with irrational input output mappings. A candidate method is modal truncation. When applying this method, we completely cut poles from the infinite partial fraction expansion of $H_{orig}$- hence, we keep only the first dominant $2N$ poles:
\begin{equation}
H_{mod}(s) = \displaystyle \sum_{k=1}^{N} \frac{4s}{s^2+c_dI\alpha_k^4s+EI\alpha_k^4}.
\end{equation}
Take $N=16$ and notice that the first peaks are accurately approximated. As for the Loewner method, the accuracy is lost outside the interest rang.

%\begin{figure}[h] \label{fig7} \vspace{-3mm}
%	\begin{center}
%		\includegraphics[scale=0.3]{plot_beam8.eps}
%		\vspace{-3mm}
%		\caption{Original and Loewner model - poles and zeros comparison}
%	\end{center} \vspace{-6mm}
%\end{figure}

A comparison between the Loewner method (applied to data coming from the original transfer function) and the modal truncation method is made by comparing the approximation error of the two methods inside the target range of frequencies. We notice that the error coming from the Loewner method is much lower than the one coming from modal (at least seven orders of magnitude), as we can see in Fig.\;7.

\begin{figure}[h] \label{fig7} \vspace{-1mm}
	\begin{center}
		\includegraphics[scale=0.28]{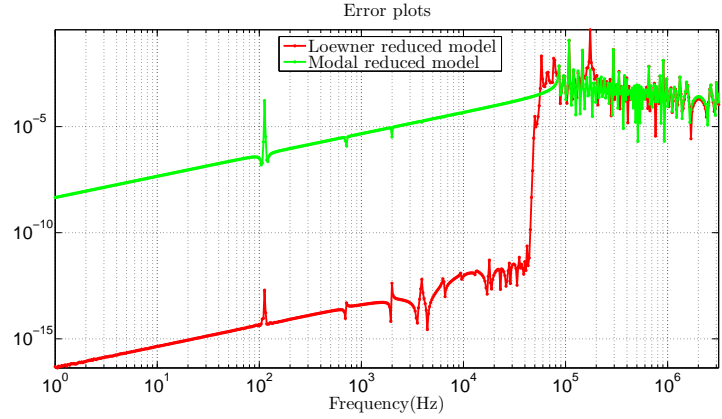}
		\vspace{-3mm}
		\caption{The error in frequency domain - Loewner vs Modal}
	\end{center} \vspace{-6mm}
\end{figure}

In the following experiments, we are going to sample the transfer function corresponding to the finite element model constructed as described in section 3. Clearly, by discretizing the spacial variable, we introduce an approximation that will considerably decrease the quality of approximation (even for large number of elements). 

Again, the first step is to take samples of $H_{FEM}(s)$ in the same 'interest' range as before. For a FE discretization with $N=1000$ elements, take as before 400 samples of the transfer function. The first major difference when using the approximated transfer function $H_{FEM}$ is the decay of singular values which is not as steep as when using samples of the original transfer function. For reduction order 32, the smallest singular value is around $10^{-6}$.

\begin{figure}[h] \label{fig8} \vspace{-3mm}
	\begin{center}
		\includegraphics[scale=0.24]{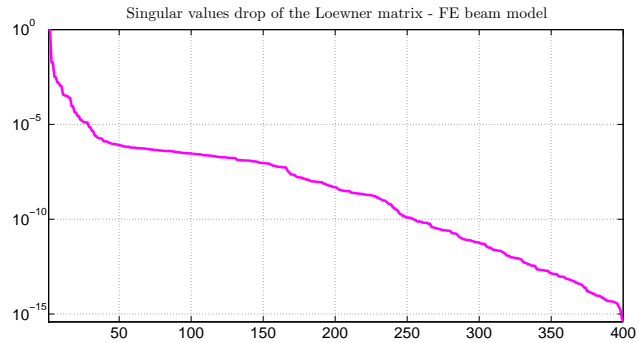}
		\vspace{-3mm}
		\caption{Singular value decay - FEM}
	\end{center} \vspace{-6mm}
\end{figure}

%The finite element approximation is comparable to the original transfer function (at least in the 'eye-ball' norm). 
%
%\begin{center}
%\includegraphics[scale=0.34]{plot_beam12.eps}
%\end{center}

Also, the poles and the zeros, are approximated well (only the four least dominant poles are completely mismatched).

\begin{figure}[h] \label{fig9} \vspace{-3mm}
	\begin{center}
		\includegraphics[scale=0.28]{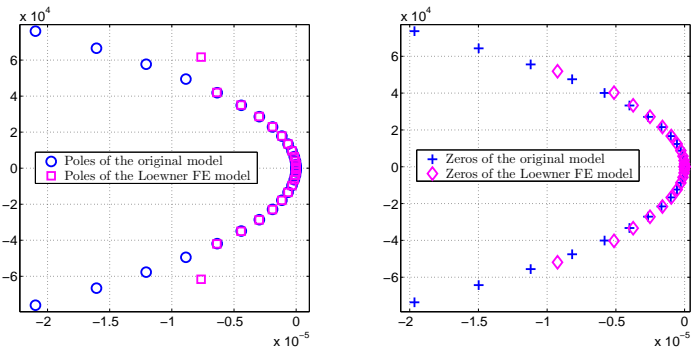}
		\vspace{-3mm}
		\caption{Poles and zeros - original and FE model}
	\end{center} \vspace{-6mm}
\end{figure}

The approximation errors are also compared. Note that the approximation quality is considerably higher when applying the Loewner framework to the data coming from the original transfer function (see Fig.\;10) as compared to the data coming  from the finite element discretization.

\begin{figure}[h] \label{fig10} \vspace{-3mm}
	\begin{center}
		\includegraphics[scale=0.28]{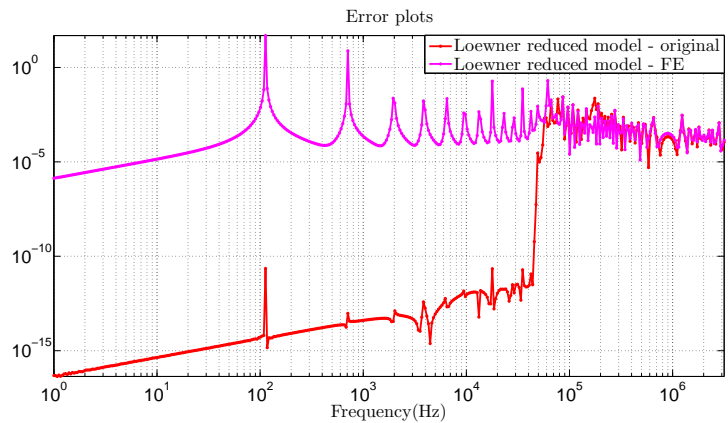}
		\vspace{-3mm}
		\caption{Error analysis}
	\end{center} \vspace{-6mm}
\end{figure}

%For the next experiment, we vary the damping parameter $c_d$ and in the same time keep the other structural parameters fixed ($E,I, L$). The next plot depicts frequency responses for various values of $c_d$. Note that, as $c_d$ increases, more and more peaks are suppressed.\\
%\noindent
%By varying the damping parameter within the range $[5 \times 10^{-4},10^2]$ it follows that the transfer functions are more or less the same in the frequency 'interest' range (at least in the 'eyeball' norm) - moreover, this behavior is noticed up to $\omega = 10^6$ Hz. 
%
%\begin{center}
%\includegraphics[scale=0.34]{plot_beam15.eps}
%\end{center}
%
%\noindent
%If we take higher and higher vallues for $c_d$, it follows that eventually more and more dominant peaks are going to be suppressed.
%
%\begin{center}
%\includegraphics[scale=0.34]{plot_beam14.eps}
%\end{center}

The Loewner framework has been recently generalized to parametric linear time invariant systems - see \cite{ail12}. In parametric model reduction, the aim is to preserve the parameters as symbolic variables in the reduced models.

After collecting frequency response measurements for appropriate
ranges of frequencies and parameter values, we use a generalization of the Loewner framework to the two variable case to construct models which are reduced both with respect to frequency and to the parameter. The Loewner matrix approach provides a trade-off between accuracy and complexity not only with respect to the frequency variable but also with respect to the parameters involved.

In this particular case, the parameter is the damping coefficient $c_d$. Our objective is to come up with a parametrized linear model that is able to faithfully reproduce the original transfer function on a particular range of frequencies as well as on a target  parameter range. As before we take measurements and we decide what order should we choose for the reduced model by inspecting the singular value decay. The frequency 'interest' range was chosen as for the other experiments, while the parameter 'target range' was taken to be $[10^{-5},10^{7}]$. For 30 values of $c_d$ in this range, we evaluate the transfer function $H_{orig}$ at 100 sampling points.

\begin{figure}[h] \label{fig11} \vspace{-3mm}
	\begin{center}
		\includegraphics[scale=0.28]{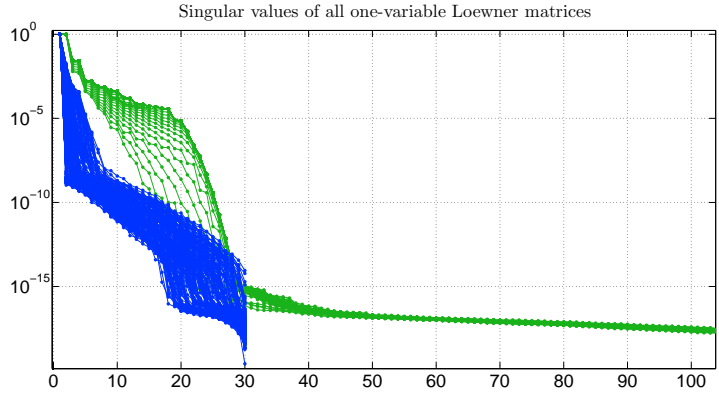}
		\vspace{-3mm}
		\caption{Singular values decay}
	\end{center} \vspace{-6mm}
\end{figure}

The plot above depicts two types of singular values. On one hand, we form Loewner matrices by using measurements for varying frequency and constant parameter ( marked with green). On the other hand, we form Loewner matrices by using measurements for varying parameter and constant frequency ( marked with blue).

By investigating the drop in the SV plot, we decide to use reduction orders $r_s = 30$ and $r_p = 14$ for building the two dimensional Loewner matrix. Then construct a reduced linear parametric model which is sampled on the same frequency and parameter range as before. When comparing to the original samples, the overall result is satisfactory:

\begin{figure}[h] \label{fig12} \vspace{-3mm}
	\begin{center}
		\includegraphics[scale=0.28]{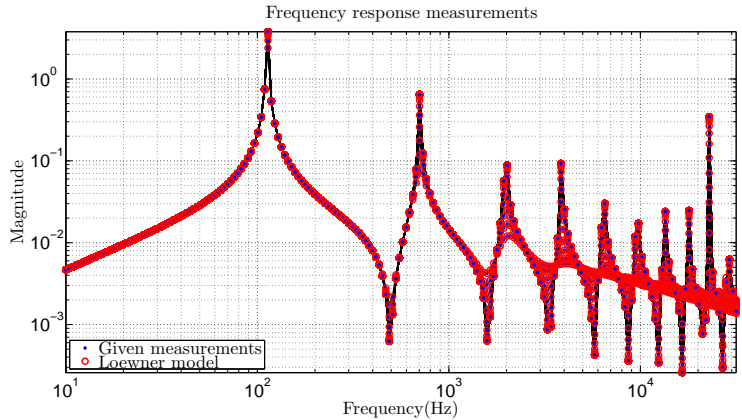}
		\vspace{-3mm}
		\caption{Comparison of the original measurements with the fitted Loewner model}
	\end{center} \vspace{-6mm}
\end{figure}

\section{Conclusion}

The aim of this study was to study the applicability of the Loewner framework to approximating a particular class of systems: 
infinite-dimensional systems. It has been shown from the various experiments performed that, indeed, the Loewner method is able to yield better approximation than other methods such as modal truncation. An advantage of the Loewner method is that it requires only input-output measurements to construct reduced order models of the data. Moreover, the use of the Loewner pencil introduces a trade-off between accuracy and complexity.

The experiments performed showed that, in terms of the accuracy of the approximation, it is better to use the framework for measurements coming directly from the original transfer function. By introducing the semi-discretization step that transforms the PDE into a series of ODE's, additional error is added that drastically decreases the approximation quality. Of course, in practice, the above mentioned step can seldom be avoided. Additionally, one could try to increase the number of elements in the FE model in order that the derivatives would be better approximated by the divided difference formulas.  

As an overall observation, the Loewner framework has been successfully applied to the new class of systems discussed in this work. As for further challenges and research interests, one would take into consideration topics such as optimal choice of sampling points or recovering from noise corrupted measurements.

\section{Appendix}

\noindent
We show the exact derivation of the finite element model for a very course discretization of the beam. For example, take $N=6$. The state equations are written as follows ($k \in \{2,3,4,5,6\}$): 
\begin{align*}
& {\color{red}(k=2) \rightarrow} \  \ddot{w}_2+EI \frac{w_{4}-4w_{3}+6w_2}{h^4} +c_dI \frac{\dot{w}_{4}-4\dot{w}_{3}+6\dot{w}_2}{h^4} = 0 \\
& {\color{red}(k=3) \rightarrow} \  \ddot{w}_3+EI \frac{w_{5}-4w_{4}+6w_3-4w_2}{h^4} +c_dI \frac{\dot{w}_{5}-4\dot{w}_{4}+6\dot{w}_3-4\dot{w}_2}{h^4} = 0 \\
& {\color{red}(k=4) \rightarrow} \
\ddot{w}_4+EI \frac{w_{6}-4w_{5}+6w_4-4w_{3}+w_{2}}{h^4} +c_dI \frac{\dot{w}_{6}-4\dot{w}_{5}+6\dot{w}_4-4\dot{w}_{3}+\dot{w}_{2}}{h^4} = 0 \\
& {\color{red}(k=5) \rightarrow} \
\ddot{w}_{5}+EI \frac{-2w_{6}+5w_{5}-4w_{4}+w_{3}}{h^4} +c_dI \frac{-2\dot{w}_{6}+5\dot{w}_{5}-4\dot{w}_{4}+\dot{w}_{3}}{h^4} = -\frac{1}{h} u \\
& {\color{red}(k=6) \rightarrow} \ \ddot{w}_{6}+EI \frac{w_{6}-2w_{5}+w_{4}}{h^4} +c_dI \frac{\dot{w}_{6}-2\dot{w}_{5}+\dot{w}_{4}}{h^4} = \frac{2}{h}u 
\end{align*}

$$
\underbrace{\left(\begin{array}{ccccc} 1 & 0 & 0 & 0 & 0\\ 0 & 1 & 0 & 0 & 0\\ 0 & 0 & 1 & 0 & 0\\ 0 & 0 & 0 & 1 & 0\\ 0 & 0 & 0 & 0 & 1 \end{array}\right)}_{{\color{blue}M}}
 \underbrace{\left(\begin{array}{c} \ddot{w}_2 \\ \ddot{w}_3 \\ \ddot{w}_4 \\ \ddot{w}_5 \\ \ddot{w}_6   \end{array}\right)}_{\ddot{v}} +
+ \underbrace{\frac{c_d I}{h^4} \left(\begin{array}{ccccc} 6 & -4 & 1 & 0 & 0\\ -4 & 6 & -4 & 1 & 0\\ 1 & -4 & 6 & -4 & 1\\ 0 & 1 & -4 & 5 & - 2\\ 0 & 0 & 1 & -2 & 1 \end{array}\right)}_{{\color{blue}J}} \underbrace{\left(\begin{array}{c} \dot{w}_2 \\ \dot{w}_3 \\ \dot{w}_4 \\ \dot{w}_5 \\ \dot{w}_6  \end{array}\right)}_{\dot{v}}+
 $$
$$ + \underbrace{\frac{E I}{h^4} \left(\begin{array}{ccccc} 6 & -4 & 1 & 0 & 0\\ -4 & 6 & -4 & 1 & 0\\ 1 & -4 & 6 & -4 & 1\\ 0 & 1 & -4 & 5 & - 2\\ 0 & 0 & 1 & -2 & 1 \end{array}\right)}_{{\color{blue}K}} \underbrace{\left(\begin{array}{c} w_2 \\ w_3 \\ w_4 \\ w_5\\w_6 \end{array}\right)}_{v} = \underbrace{\frac{1}{h}\left(\begin{array}{c}  0\\ 0\\ 0\\ -1 \\ 2  \end{array}\right)}_{{\color{blue}f}} u.
$$
Next, by applying the Laplace transform of the above time-domain collection of ODE's, we obtain:
$$
s^2MV(s)+sJV(s)+KV(s) = fU(s) \Rightarrow V(s) = (s^2M+sJ+K)^{-1}fU(s).
$$ 
We derive the output equation by using the surrogate variable $z$ as follows:
$$
EI z + c_d I \dot{z} =  \underbrace{EI \left(\begin{array}{ccccc} 0 & 0 & 0 & -2 & 3 \end{array}\right)
}_{{\color{blue}c_1}} \underbrace{\left(\begin{array}{c} w_2 \\ w_3 \\ w_4 \\ w_5\\w_6 \end{array}\right)}_{v} + \underbrace{c_d I \left(\begin{array}{ccccc} 0 & 0 & 0 & -2 & 3 \end{array}\right)
}_{{\color{blue}c_2}} \underbrace{\left(\begin{array}{c} \dot{w}_2 \\ \dot{w}_3 \\ \dot{w}_4 \\ \dot{w}_5 \\ \dot{w}_6 \end{array}\right)}_{\dot{v}} + \underbrace{2h^3}_{{\color{blue}d}} u.
$$
Now, by taking the Laplace transform of the above differential equation in $z$, we obtain that: 
$$
EI Z(s) + c_dIsZ(s) = c_1 V(s) + c_2 s V(s)+ dU(s) \Rightarrow Z(s) = \frac{(c_1+sc_2)V(s)}{EI+sc_dI}+\frac{d}{EI+sc_dI}U(s).
$$
Since $z = \dot{w}_{N+2}$ and $y = w_{N+2}$ it follows that $z = \dot{y}$ and by taking again Laplace transform we have: $sZ(s) = Y(s) \Rightarrow Z(s) = \frac{Y(s)}{s}$. Substituting $Z(s)$ in the above equation we get:
$$
Y(s) = \frac{s(c_1+c_2s)V(s)}{EI+c_dS}+\frac{ds}{EI+c_dS}U(s).
$$
Substituting $V(s)$ in the above equation, we form the transfer function of the FEM as follows:
$$
H_{FEM}(s) = \frac{Y(s)}{U(s)} = \frac{s}{EI+sc_dI} \Big{[} (c_1+c_2s) (s^2M+sJ+K)^{-1}f+d \Big{]}.
$$ 
By introducing the new augmented state vector $x = \left(\begin{array}{c} v \\ \dot{v}  \end{array}\right) \ \in \mathbb{R}^{2(n-1)}$, we rewrite the dynamics of the second degree differential system above as a first degree dynamical system in $x$:
$$
G \dot{x} = A x + B u , \ \ EIz+c_dIy = C x + D u,
$$
where:
$$
G = \left(\begin{array}{cc} I_{n-1} & 0 \\ 0 & M   \end{array}\right), \ A = \left(\begin{array}{cc} 0 & I_{n-1} \\ -K & -J  \end{array}\right), \ B = \left(\begin{array}{c} 0 \\ f  \end{array}\right), \ C = \left(\begin{array}{cc} c_1 & c_2   \end{array}\right), \ D = d,
$$
and $I_{n-1}$ is the identity matrix of dimension $n-1$. Taking the Laplace transform we get:
$$
X(s) = (sG-A)^{-1}B U(s) \ \ \text{and} \ \ EIZ(s)+c_dIY(s) = CX(s)+DU(s).
$$
Again, using that: $sZ(s) = Y(s)$, it follows that:
$$
H_{FEM}(s) = \frac{Y(s)}{U(s)} = \frac{s}{EI+sc_dI} [C(sG-A)^{-1}B+D].
$$

\end{document}